\documentclass[11pt]{amsart}
\usepackage{amsfonts}
\usepackage{amsmath}
\usepackage{amsthm}
\usepackage{bm,bbm}
\usepackage[dvips]{graphicx}
\usepackage{caption}
\usepackage{color}
\usepackage{multicol}

\usepackage{bigints}

\usepackage{mathtools}
\usepackage{commath}

\numberwithin{equation}{section}

\allowdisplaybreaks

\theoremstyle{plain}
\newtheorem{theorem}{Theorem}[section]
\newtheorem{lemma}[theorem]{Lemma}
\newtheorem{proposition}[theorem]{Proposition}

\theoremstyle{definition}

\def\rl{r_\la(s_0)}
\def\rlw{r_\la(w_\la)}
\def\bet{\begin{theorem}}
\def\ent{\end{theorem}}
\def\pa{\partial}
\def\om{\omega}

\def\z{\zeta}
\def\t{\tau}
\def\ti{\times}

\def\k{\kappa}

\def\ff{\infty}
\def\tff{\uparrow\ff}

\def\la{\leftarrow}
\def\f{\frac}
\def\d{\dif}
\def\de{\delta}
\def\ms{\mathsf}
\def\beqn{\begin{equation}}
\def\beqn*{$$}
\def\eeqn{\end{equation}}

\newcommand{\su}{\subseteq}
\newcommand{\sm}{\setminus}

\newcommand{\bx}{{\bf x}}

\def\mc{\mathcal}
\def\P{\mathbb{P}}
\def\PP{\mc{P}}
\def\QQ{\mc Q}
\def\co{\colon}

\def\E{\mathbb{E}}
\def\Pn{\mathcal P}

\newcommand{\reals}{{\mathbb R}}

\newcommand{\R}{\reals}

\newcommand{\X}{{\mathcal{X}}}

\newcommand{\I}{\mathcal I}

\newcommand{\one}{{\mathbbm 1}}

\newcommand{\remove}[1]{}

\newcommand{\Z}{\mathbb Z}

\def\enp{\end{proof}}
\def\bel{\begin{lemma}}
\def\bep{\begin{proof}}
\def\enl{\end{lemma}}
\newcommand{\M}{\mathcal M}

	\def\hek{\hat{\eta}_{\la}}
	\def\hzk{\hat{\zeta}_{\la}}
	\def\ek{\eta_{\la}}
	\def\nla{N^{(1)}_\la}
	\def\nlb{N^{(2)}_\la}
	
	\def\dll{d_{\ell, \la}}
	\def\ekl{\eta^{(m)}_{ \la}}
	\def\hekl{\hat{\eta}^{(m)}_{ \la}}
	\def\hzkl{\hat{\zeta}^{(m)}_{ \la}}
	\def\heol{\hat{\eta}^{(1)}_{ \la}}
	\def\hzol{\hat{\zeta}^{(1)}_{ \la}}
	
	\def\zk{\zeta_{ \la}}
	\def\zkl{\zeta^{(m)}_{ \la}}
	\def\been{\begin{enumerate}}
	\def\enen{\end{enumerate}}
\def\im{\item}
\def\im{\item}

\renewcommand{\H}{\mathbb H}
\renewcommand{\la}{\lambda}
	\newcommand{\bepr}{\begin{proposition}}
	\newcommand{\enpr}{\end{proposition}}
	\newcommand{\dist}{\mathsf{dist}}
\def\LL{\mc L}

\definecolor{OliveGreen}{cmyk}{0.64,0,0.95,0.40}

\begin{document}

\bibliographystyle{abbrv}

\title[Large deviations for hyperbolic $k$-nearest neighbor balls]
{Large deviations  for the volume \\ of hyperbolic $k$-nearest neighbor balls}

\author{Christian Hirsch}
\author{Moritz Otto}
\address[Christian Hirsch, Moritz Otto]{Department of Mathematics, Aarhus University, Ny Munkegade 118, 8000 Aarhus C, Denmark}
\email{hirsch@math.au.dk, otto@math.au.dk}
\address[Christian Hirsch]{DIGIT Center, Aarhus University, Finlandsgade 22, 8200 Aarhus N, Denmark}

\author{Takashi Owada}
\address[Takashi Owada]{Department of Statistics\\
Purdue University \\
West Lafayette, 47907, USA}
\email{owada@purdue.edu}

\author{Christoph Th\"ale}
\address[Christoph Th\"ale]{Faculty of Mathematics\\ Ruhr University Bochum \\   44780, Bochum, Germany.}
\email{christoph.thaele@rub.de}

\subjclass[2010]{Primary 60F10. Secondary 51M10, 52A55, 60D05, 60G55.}
\keywords{hyperbolic space, large deviation principle, nearest neighbor balls, Poisson point process, stochastic geometry.}

\begin{abstract}
	We prove a large deviation principle for the point process of large Poisson $k$-nearest neighbor balls in hyperbolic space. More precisely, we consider a stationary Poisson point process of unit intensity in a growing sampling window in hyperbolic space. We further take a growing sequence of thresholds such that there is a diverging expected number of Poisson points whose $k$-nearest neighbor ball has a volume exceeding this threshold. Then, the point process of exceedances satisfies a large deviation principle whose rate function is described in terms of a relative entropy. The proof relies on a fine coarse-graining technique such that inside the resulting blocks the exceedances are approximated by independent Poisson point processes.
\end{abstract}
\maketitle

%
%
\section{Introduction}
\label{sec:intro}

%
%
Large deviations theory is one of the most classical subfields of probability theory with a wide range of applications in information theory, statistical physics, and rare-event simulation \cite{dz98}. While the most refined results are available in the study of sequences of random variables and of time-dependent stochastic processes, the understanding of large deviations in spatial random systems is still in its infancy. Even the most basic statistics such as the edge counts in a random geometric graph can give rise to surprising and highly non-trivial condensation effects on the level of large deviations \cite{harel}.

%
%
One of the earliest achievements in this context is \cite{georgii2}, which  proves a large deviation principle (LDP) for statistics of a homogeneous Poisson point process in $\R^d$ under rather restrictive boundedness and locality assumptions. Later, these conditions could be relaxed substantially so as to allow for statistics with bounded exponential moments and stabilizing score functions  \cite{yukLDP,yukLDP2}. While the main focus there was on scalar and measure-valued LDPs in the thermodynamic regime, very recently also geometric statistics in the sparse and dense regimes of the Poisson point process were considered \cite{kang,gtf}.

%
%
All of the aforementioned results have in common that they deal with random geometric systems  in Euclidean space. However, during the last years the hyperbolic space has received substantial attention in the context of complex networks \cite{fount}. Moreover, there has been vigorous activity to understand the asymptotic behavior of geometric functionals of  random set systems in hyperbolic space as well. While there has been substantial progress for central and non-central limit theorems and Poisson approximation results \cite{BesauRosenThaele,BesauThaele,BetkenHugThaele,flammant,FountoulakisYukich,GKTBetaStar,HeroldHugThaele,KRTlambda,otto:thale:2022,subtree}, large deviation principles have not been considered so far. By investigating the volume of $k$-nearest neighbor (kNN) balls in hyperbolic space in the present work, we provide a first step in this direction and complement the recent findings \cite{otto:thale:2022} about their extremal behavior.

%
%
The general proof strategy is to extend and to develop further the ideas of \cite{kang} on a coarse graining scheme to introduce a blocked point process. Inside each block, a Poisson approximation theorem in the spirit of \cite{otto:thale:2022} is used to replace the original functional by a Poisson point process for which the LDP is given by Sanov's theorem. However, the intricate geometry of the hyperbolic space makes it far more challenging to implement the blocking argument in comparison to the Euclidean setting considered in \cite{kang}. For instance, in the Euclidean setting the kissing number is finite so that in a box of side length of order $r > 0$ only a uniformly bounded number of points with a nearest neighbor radius exceeding $r$ can be placed. In contrast, in hyperbolic space, this number grows exponentially in $r$; see \cite{DostertKolpakov}.  Moreover, in the half-space model of  hyperbolic space, the Poisson intensity is inhomogeneous in the vertical direction when considered in Euclidean terms, while in the $\R^d$-setting there is homogeneity in all directions. We deal with these problems by deriving more refined exponential moment bounds and analyzing the point configuration in vertical layers that individually can be considered approximately homogeneous.

%
%
We believe that the techniques developed in the present article open the door to the investigation of further LDPs in hyperbolic space. In particular, we think of extending the results for the Euclidean component counts from \cite{gtf} in the sparse regime.  Here, we note that while the study of component counts is restricted to trees in the dense regime \cite{subtree} such constraints are no longer present in the sparse regime.

%
%
The rest of the manuscript is organized as follows. In Section \ref{sec:mod}, we properly introduce the considered model and state the LDP for the volumes of large kNN balls. Next, in Section \ref{sec:out}, we give a proof outline, where we reduce the assertion to two key auxiliary results, namely Propositions \ref{pr:ee1} and \ref{pr:ee2}. These results are proven separately in Sections \ref{sec:ee1} and \ref{sec:ee2}, respectively.

%
%
\section{Model and main result}
\label{sec:mod}

%
%
By the hyperbolic space $\H^d$ of dimension $d\geq 2$ we mean the unique simply connected, $d$-dimensional Riemannian manifold of constant sectional curvature $-1$, cf.\ \cite{flavors,ratcl}. There are several models to represent $\H^d$ in the $d$-dimensional Euclidean space $\R^d$. To carry out our computations, we will work in this paper with the so-called half-space model, but we emphasize that all results we derive are actually model independent; for background material on half-space model of hyperbolic space we refer to \cite[Chapter 4.6]{ratcl}.  In the half-space model, we identify $\H^d$ with the product space $\H^d = \R^{d - 1} \ti (0, \ff)$. The Riemannian metric is then determined by
$$
\d s^2 = \f{\d x_1^2 + \ldots + \d x_{d - 1}^2 + \d y^2}{y^2},\qquad(x_1,\ldots,x_{d-1})\in\R^{d-1},\,y\in(0,\infty),
$$
and we denote by $\dist_{\ms{hyp}}(z_1,z_2)$ the hyperbolic distance of two points $z_1=(x_1,y_1),z_2=(x_2,y_2)\in\H^d$, which in terms of the Euclidean distance $\dist_{\ms{euc}}(z_1,z_2)$ between $z_1$ and $z_2$ is given by
$$
\dist_{\ms{hyp}}(z_1,z_2) = {\rm arcosh}\Big(1 + {\dist_{\ms{euc}}(z_1,z_2)^2\over 2y_1y_2}\Big),
$$
see \cite[Theorem 4.6.1]{ratcl}.
According to \cite[Theorem 4.6.7]{ratcl}, the hyperbolic volume measure $V_{\ms{hyp}}$ on $\H^d$ has density
\begin{align}
	\label{eq:vol}
	{\d x_1 \cdots \d x_{d - 1} \d y\over y^d}
\end{align}
with respect to the Lebesgue measure on $\R^{d-1}\times(0,\infty)$ and we will use the notation $$|B|_{\ms{hyp}}:=V_{\ms{hyp}}(B)=\int_{B}y^{-d}\d x_1 \cdots \d x_{d - 1} \d y$$ for the hyperbolic volume of a measurable set $B\subseteq\H^d$. In contrast, we shall write $|\,\cdot\,|_{\ms{leb}}$ for the Lebesgue measure of the appropriate dimension, which will always be clear from the context. In addition, we denote for $r>0$ and $x\in\H^d$ by $B_{r}(x):=\{z\in\H^d:\dist_{\ms{hyp}}(x,z)\leq r\}$ the hyperbolic ball of radius $r$ centered at $x$. Its hyperbolic volume satisfies
\begin{equation}\label{eq:VolBall}
|B_r(x)|_{\ms{hyp}}=\beta_d\int_0^r\sinh^{d-1}u\,\d u,
\end{equation}
independently of $x$, with $\beta_d:=2\pi^{d/2}/\Gamma({d\over 2})$ being the surface content of the $(d-1)$-dimensional Euclidean unit sphere; see \cite[page 79]{ratcl}. In particular, $|B_r(x)|_{\ms{hyp}}$ grows like a constant multiple of $e^{r(d-1)}$, as $r\uparrow\infty$.

Henceforth, we let $\PP = \sum_{i \ge 1} \de_{X_i}$ be the random counting measure distributed as a Poisson point process with the hyperbolic volume measure as its intensity measure. We note that such $\PP$ is stationary in the sense that its distribution is invariant with respect to the full group of hyperbolic isometries. In particular, $\PP(B)$ is the number of points of $\PP$ falling into a measurable set $B\subseteq\H^d$.
In this paper, we write  $\PP_B$ or sometimes $\PP\cap B$ for the restriction of the measure $\PP$ to $B$. In what follows, we shall restrict $\PP$ to the family of sampling windows
$$W_\la := [0, 1]^{d - 1} \ti [e^{-\la}, \ff),\qquad\la>0,$$
whose hyperbolic volume is given by $|W_\la|_{\ms{hyp}} = \int_{e^{-\la}}^\infty y^{-d}\,\d y= {1\over d-1}e^{(d - 1)\la }$. 

%
%
For $k\geq 1$ we study the asymptotics of large $k$-nearest neighbor (kNN) balls centered in $W_\la$. To make this precise, we define the point process 
$$
\xi_{k, \la} := \sum_{x \in \PP_{W_\la}}\delta_{|B_{R_{ k}(x)}(x)|_{\ms{hyp}} - v_\la},
$$
where
$$R_k(x):= \inf\{r \ge 0 \co \PP(B_r(x)) \ge k + 1\}$$
and $(v_\la)_{\la>0}$ is a threshold sequence satisfying 
$$v_\la -(d - 1)\la - (k - 1)\log \la \to -\ff, \qquad \la \tff.$$
In particular, the expected number of exceedances in the window $W_\la$ is of order $u_\la := |W_\la|_{\ms{hyp}}e^{ - v_\la} v_\la^{k - 1}$. 

%
%
Our main result is a large deviation principle (LDP) for the volumes of large kNN balls. We recall that a family of random variables $(X_\lambda)_{\lambda>0}$, defined on some probability space $(\Omega,\mathcal{F},\P)$ and taking values in a Polish space $\X$, satisfies an LDP with speed $s_\lambda\uparrow\infty$ and rate function $\I:\X\to[0,\infty]$, provided that $\I$ is lower semicontinuous and if for each measurable set $B\subseteq\X$ one has that
\begin{align*}
-\inf_{x\in B^\circ}\I(x) \leq\liminf_{\la\uparrow\infty}s_\la^{-1}\log\P(X_\la\in B)\leq\limsup_{\la\uparrow\infty}s_\la^{-1}\log\P(X_\la\in B)\leq-\inf_{x\in\bar B}\I(x),
\end{align*}
where $B^\circ$ and $\bar B$ stand for the interior and the closure of $B$, respectively. To present our main result, we fix $s_0\in\R$ and introduce the space $E_0 := [s_0, \ff)$ as well as the measure $\t_k$ on $E_0$ which has Lebesgue density $e^{-u}/(k - 1)!$, $u>s_0$. By $\M(E_0)$ we denote the space of finite measures on $E_0$, which supplied with the weak topology becomes a Polish space \cite[Proposition A2.5.III]{pp1}. For $\rho\in\M(E_0)$ we let 
\begin{equation}  \label{e:def.rela.entropy}
H(\rho|\tau_k) := \begin{cases}
\int_{E_0}\log \f{\dif \rho}{\dif \tau_k}(\bx)\rho(\dif \bx) -\rho(E_0) +\tau_k(E_0) & \text{ if } \rho\ll \tau_k, \\
\infty & \text{ otherwise}
\end{cases}
\end{equation}
be the relative entropy entropy of $\rho$ with respect to $\t_k$, where $\rho\ll \tau_k$ indicates that $\rho$ is absolutely continuous with respect to $\t_k$ and $\f{\dif \rho}{\dif \tau_k}$ denotes the corresponding Radon-Nikodym derivative; see \cite[Equation (2.10)]{georgii2}. We are now prepared to present the main result of this paper.

\bet[LDP for kNN balls]
\label{thm:knn}
Let $k \ge 1$. Then, the family of random measures $\xi_{k, \la}/u_\la$ satisfies an LDP on $\M(E_0)$ with speed $u_\la$ and rate function $H(\,\cdot\, |\t_k)$.
\ent

As highlighted in Section \ref{sec:intro}, Theorem \ref{thm:knn} can be seen as the hyperbolic counterpart of \cite[Theorem 2.1]{kang}, which concerns large deviations of the empirical measure of recentered and rescaled kNN balls in Euclidean space. On a formal level, in the hyperbolic setting, we found it more convenient to consider a scaling where the expected number of Poisson points in the window grows exponentially in $\la$, whereas in \cite{kang} the corresponding scaling is linear in parameter $n$. After this rescaling, our condition on the growth of $v_\la$ corresponds precisely to the growth of $a_n$ in \cite[Equation (2.1)]{kang}. Moreover, in the Euclidean setting the dense scaling in \cite{kang} can be equivalently transformed to the regime of growing windows considered here. However, we stress that the regimes of dense points and growing windows are no longer equivalent in hyperbolic setting, and only the growing-window asymptotics reflects the negative curvature effects from the hyperbolic space.

While the previous paragraph illustrates that there is a formal correspondence between Theorem \ref{thm:knn} and the Euclidean analog \cite[Theorem 2.1]{kang}, the hyperbolic geometry creates substantial complications, which require us to develop novel methodological tools. First, while both Theorem \ref{thm:knn} and \cite[Theorem 2.1]{kang} rely on a Poisson approximation result for large kNN balls, such a result is substantially harder to establish in hyperbolic geometry. We therefore adapt the arguments of a very recent work \cite{otto:thale:2022}. Moreover, when establishing a central uniform integrability property, \cite{kang} relies crucially on a kissing-number argument. More precisely, while in Euclidean space, the number of non-intersecting balls of radius $r$ that can be put into a box of side length of order $r$ is uniformly bounded independently of $r$, such a property does not hold in hyperbolic space. Therefore, we need to develop delicate exponential moment bound that ensure uniform integrability despite a potentially unbounded number of balls.

To make the presentation more accessible, we provide a brief summary of the notation used throughout this paper, most terms will formally be defined in the forthcoming sections:
\begin{multicols}{2}
	\noindent- $\#(\,\cdot\,)$: cardinality of a set\\
	- $\one\{\,\cdot\,\}$: indicator function \\
	- $o(\,\cdot\,),O(\,\cdot\,)$: Landau symbols\\	
	- $(\Omega,\mathcal{F},\P)$: underlying probability space\\
	- $\E[\,\cdot\,]$: expectation(integration) wrt.\ $\P$\\
	- $\mathcal{L}(\,\cdot\,)$: law of a random element\\
	- $\delta_{(\,\cdot\,)}$: Dirac measure\\
	- $\M(\,\cdot\,)$: space of locally finite measures \\
	- $\H^d$: $d$-dimensional hyperbolic space,\\
	\phantom{x}\qquad identified with $\R^{d-1}\times(0,\infty)$\\
	- $H(\,\cdot\,|\tau)$: relative entropy of a measure wrt.\ $\tau$\\	
	- $\omega_A$: restriction of a measure to a set $A$\\
	- $B_r(x)$:hyperbolic ball of radius $r$ centered at $x$\\
	- $|\,\cdot\,|_{\ms{hyp}}$: hyperbolic volume measure\\ \phantom{x}\hspace{1.3cm} with differential $\d V_{\ms{hyp}}$\\
	- $|\,\cdot\,|_{\ms{leb}}$: Lebesgue measure\\
	- $\dist_{\ms{hyp}}(\,\cdot\,,\,\cdot\,)$: hyperbolic distance\\
	- $\dist_{\ms{euc}}(\,\cdot\,,\,\cdot\,)$: Euclidean distance\\
	- $\dist_{\ms{TV}}(\,\cdot\,,\,\cdot\,)$: total variation distance\\
	- $\dist_{\ms{KR}}(\,\cdot\,,\,\cdot\,)$: Kantorovich-Rubinstein distance
\end{multicols}

In addition, simple point processes $\omega=\sum_{k\geq 1}\delta_{z_k}$ with $z_k\in\H^d$ will be identified with locally finite point configurations $\{z_k:k\geq 1\}$. Moreover, by abuse of notation, we write $z \in \omega$ if $z$ is an atom of the counting measure $\omega$. The intensity measure of $\omega$ is denoted by $\E[\omega]$.

In this paper we will denote by $C>0$ a generic constant whose value might change from occurrence to occurrence. If several such constants are needed at the same time, we denote them by $C_0, C_1,\ldots$

Finally, since the parameters $k$ and $s_0$ are fixed in the rest of the manuscript, we omit them from the notation.

%
%
\section{Proof outline}
\label{sec:out}
%
%
To prove Theorem \ref{thm:knn}, we partition $[0, 1]^{d - 1}$ into $u_\la$ congruent blocks $S_1, \dots, S_{u_\la}$ and set $Q_m := S_m \ti [e^{-\la}, \ff)$. We note that $u_\la$ is not necessarily an integer. However, to keep notation simple we do not write rounding symbols.
The key step in the proof of Theorem \ref{thm:knn} is to relate the empirical measure $\xi_{ \la}/u_\la$ with the following separated point process: 
$$
\ek := \sum_{m \le u_\la}\ekl,\qquad\text{where}\qquad\ekl:=\sum_{x \in \PP_{Q^-_m}}\one\{f(x, \PP_{Q_m}) > s_0\}\delta_{f(x, \PP_{Q_m})}.
$$
Here, for a locally finite counting measure $\omega$ on $\H^d$ we put $f(x, \om) := |B_{R_k(x, \om)}|_{\ms{hyp}} - v_\la$ with $R_k(x,\omega):=\inf\{r\geq 0:\omega(B_r(x))\geq k+1\}$ and let $Q^-_m \su Q_m$ be an `internal region' that we will define precisely in \eqref{eq:km} below. The idea behind introducing these internal regions is that we want large kNN balls  to occur independently for distinct regions. In other words, this step will remove the dependence of kNN radii exceeding the threshold $r_\la(u)$, which is  defined  such that $|B_{r_\la(u)}|_{\ms{hyp}} = u + v_\la$.

%
%
Having introduced the process $\ek$, the proof of Theorem \ref{thm:knn} can be split up into two steps regarding exponential equivalence (with respect to the total variation distance) in the sense of \cite[Definition 4.2.10]{dz98}:
\been
\im show that the family of random measures $\ek/u_\la$ is exponentially equivalent to  $\z_{ \la}/u_\la$, where for each $\lambda>0$, $\z_{ \la}$ is a Poisson point process on $E_0$ whose intensity measure has density $u_\la\tau_k$ with respect to the Lebesgue measure;
\im show that the family of  random measures $\ek/u_\la$ is exponentially equivalent to $\xi_{ \la}/u_\la$.
\enen

%
%
\bepr[Exponential equivalence of $\ek$ and $\z_{ \la}$]
\label{pr:ee1}
The families of random measures $\ek/u_\la$ and $\z_{ \la}/u_\la$ are exponentially equivalent.
\enpr

%
%
\bepr[Exponential equivalence of $\ek$ and $\xi_{ \la}$]
\label{pr:ee2}
The families of random measures $\ek/u_\la$ and $\xi_{ \la}/u_\la$ are exponentially equivalent.
\enpr

%
%
The main task is to prove Propositions \ref{pr:ee1} and \ref{pr:ee2}, which are the hyperbolic analogs of \cite[Proposition 4.3]{kang} and \cite[Propositions 4.4--4.6]{kang}. Before doing so in Sections \ref{sec:ee1} and \ref{sec:ee2}, we briefly explain how to deduce Theorem \ref{thm:knn} from these two results.

\bep[Proof of Theorem \ref{thm:knn}]
First, we note that by the Poisson variant of Sanov's theorem, the family of random measures $(\z_{\la}/u_\la)_{\la>0}$ satisfies an LDP with speed $u_\la$ and rate function $H(\,\cdot\, |\t_k)$; see \cite[Theorem 6.2.10]{dz98}. Hence, taking into account Propositions \ref{pr:ee1} and \ref{pr:ee2}, the result follows from the fundamental fact in large deviations theory that exponentially equivalent families of random elements satisfy the same LDP; see \cite[Theorem 4.2.13]{dz98}.
\enp

We conclude the present overview section, with a precise definition of the internal regions $Q^-_m$. To that end, we  let $(w_\la)_{\la>0}$ be a diverging sequence with $w_\la \in o(v_\la)$. Now, for each $y \ge e^{-\la}$, we define 
 \begin{align}
	 \label{eq:km}
	 Q^-_m := \{(x, y) \in Q_m:  x \in S^-_{m ,y} \}
 \end{align}
where 
$S^-_{m, y} := \{x \in S_m \co \dist_{\ms{euc}}(x,\pa S_m) \ge x_\la(y)\} \subseteq \R^{d - 1}$ with
$x_\la(y) := ye^{\rlw }$, $\dist_{\ms{euc}}$ refers to the Euclidean distance in $\R^{d-1}$ and $\pa S_m$ stands for the boundary of the set $S_m$.
The desired properties of the internal regions $Q^-_m$ are described in the following proposition.
\bepr[Containment and volume of $Q^-_m$]
\label{pr:cvkm}
It holds that 
\been 
\im $B_{r_\la(w_\la)}(z) \su Q_m$ for every $z \in Q^-_m$.
\im $\big|Q_m \setminus Q^-_{m}\big|_{\ms{hyp}} \in o(|Q_m|_{\ms{hyp}})$.
\enen
\enpr

To prove Proposition \ref{pr:cvkm}, we rely on the following distance formula taken from \cite[Proposition 2.1]{flammant}. We also recall that we identify points $z\in\H^d$ with pairs  $z=(x,y)$ in the product space $\R^{d-1}\times(0,\infty)$.

%
%
\bel
\label{lem:dist}
Let $z_1=(x_1,y_1)\in\H^d$, $z_2=(x_2,y_2)\in\H^d$ and define $\k := \dist_{\ms{euc}}(x_1,x_2)/y_1$ and $v := y_2/y_1$. Then, 
$$\dist_{\ms{hyp}}(z_1, z_2) := \Phi(v^{-1}(\k^2 + (v + 1)^2),$$
where 
$\Phi(t) := \log(t) - \log(4) + 2\log( 1+ \sqrt{1 - 4/t}).$
In particular, by minimizing over $y_2$, 
\begin{align}
	\label{eq:lk}
	\inf_{y_2 > 0}\dist_{\ms{hyp}}(z_1, (x_2, y_2)) \ge \log(\k).
\end{align}
holds for all sufficiently large $\k > 0$.
\enl
\bep
To carry out the minimization, we set the $\tfrac \d{\d v}(v^{-1}(\k^2 + (v + 1)^2)) = 0$. This gives $1 = (\k^2 + 1)/v^2$. Hence, $v^{-1}(\k^2 + (v + 1)^2)$ is asymptotically equivalent to $2 \sqrt{\k^2 + 1}$, which in turn is equivalent to $2\k$ for large $\kappa$.
\enp

We now prove Proposition \ref{pr:cvkm}. For the proof, we note that 
\begin{equation}  \label{e:bound.r.la.w.la}
\frac{\log (v_\la)}{d-1} +C_1 \le r_\la (w_\la) \le \frac{\log (v_\la)}{d-1} +C_2, 
\end{equation}
for suitable constants $-\infty<C_1<C_2<\infty$; see Equation (4.1) in \cite{otto:thale:2022}.

\bep[Proof of  Proposition \ref{pr:cvkm}] 
We prove separately parts (i) and (ii).
\medskip

\paragraph{Part (i)} Let $z_1 = (x_1, y_1) \in Q^-_m$ and $x_2 \in  S_m^\circ$. Then, by \eqref{eq:lk},
$$\min_{y_2 > 0}\dist_{\ms{hyp}}(z_1, (x_2, y_2))\ge \log(\k)  \ge \log(x_\la(y_1)/y_1)  = r_\la(w_\la).$$
In other words, $B_{r_\la(w_\la)}(z_1) \su Q_m$.
\medskip

\paragraph{Part (ii)}
Note that $Q_m$ is a cube of side length $u_\la ^{-1/(d - 1)}$ so that  $$|S_m \sm S^-_{m, y}|_{\ms{leb}} \le 2(d - 1)u_\la^{-(d - 2)/(d - 1)}x_\la(y).$$ Hence, for $d \ge 3$, by \eqref{eq:vol} and Fubini's theorem,  $|Q_m \sm Q^-_m|_{\ms{hyp}}$ can be bounded as
\begin{align*}
 \int_{e^{-\la}}^\ff y^{-d}|S_m \sm S^-_{m, y}|_{\ms{leb}} \d y 
	&\le2(d - 1)u_\la^{-(d - 2)/(d - 1)} \int_{e^{-\la}}^\ff x_\la(y)y^{-d} \d y \\
	&\le C u_\la^{-(d - 2)/(d - 1)}v_\la^{1/(d - 1)}|W_\la|_{\ms{hyp}}^{(d - 2)/(d - 1)}.
\end{align*}
 Therefore,
$$\f{\big|Q_m \sm Q^-_m\big|_{\ms{hyp}}^{d - 1}}{|Q_m|_{\ms{hyp}}^{d - 1}} \le \f{c u_\la^{-(d - 2)}v_\la|W_\la|_{\ms{hyp}}^{d - 2}}{u_\la^{-(d - 1)}|W_\la|_{\ms{hyp}}^{d - 1}} = C e^{-v_\la/(d - 1)}v_\la^{k/(d - 1)}\to0,$$
as $\la\uparrow\infty$.

Next, consider $d = 2$. Here, the important observation is that $S^-_{m, y} = \emptyset$ if $x_\la(y) \ge u_\la^{-1}$. That is, if $y \ge y_0(\la) := u_\la^{-1}e^{-r_\la(w_\la)}$. Now,
$$\int_{y_0(\la)}^\ff y^{-2} |S_m|_{\ms{leb}}\d y = {u_\la^{-1}y_0(\la)}^{-1} = e^{r_\la(w)} \le Cv_\la \in o\big(u_\la^{-1}|W_\la|_{\ms{hyp}}\big).$$
Moreover,
\begin{align*}
	\int_{e^{-\la}}^{y_0(\la)}y^{-2}|S_m \sm S^-_{m, y}|_{\ms{leb}} \d y 
	&\le2(d - 1) \int_{e^{-\la}}^{y_0(\la)} x_\la(y)y^{-2} \d y\\
	&\le C v_\la (\log(y_0(\la)|W_\la|_{\ms{hyp}}))\\
	&\le C_1 v_\la (v_\la +C_2\log(v_\la)).
\end{align*}
Noting that  $v_\la^2 \in o\big(u_\la^{-1}|W_\la|_{\ms{hyp}}\big)$ concludes the proof.
\enp

%
%
\section{Proof of Proposition \ref{pr:ee1}}
\label{sec:ee1}

In this section, we prove Proposition \ref{pr:ee1}, that is, the exponential equivalence of the empirical measures $\ek/u_\la$ and a family of empirical measures of suitable Poisson point processes. We recall that $\ek$ is composed of its constituents $\ekl$ in the individual boxes. The first step of the proof is therefore to proceed as in \cite[Proposition 4.1]{kang} and establish a Poisson approximation result individually for each $\ekl$. We do this by showing the stronger assertion that for each $m\geq 1$ the Kantorovich-Rubinstein distance $\dist_{\ms{KR}}(\LL(\ekl), \LL(\zkl))$ between the law $\LL(\ekl)$ of $\ekl$ and that of $\zkl$ tends to zero, as $\la\uparrow\infty$.  Here, we recall that for two point processes $\omega_1,\omega_2$ on $\H^d$, $$\dist_{\ms{KR}}(\LL(\omega_1),\LL(\omega_2)):=\sup_h|\E h(\omega_1)-\E h(\omega_2)|,$$ where the supremum runs over all measurable Lipschitz-$1$ functions on the space of point processes on the space $E_0$ with respect to the total variation distance $\dist_{\ms{TV}}$. For two measures $\mu_1,\mu_2$ on $E_0$ the latter is given by $$\dist_{\ms{TV}}(\mu_1,\mu_2):=\sup_{A\subseteq E_0}|\mu_1(A)-\mu_2(A)|.$$

%
%
\bepr[Poisson approximation for separated processes]
\label{lem:sep} 
Let $m \ge 1$. Then,
$$\dist_{\ms{KR}}\big(\LL(\ekl), \LL(\zkl)\big)\to 0,\quad\text{as $\la\tff$},$$
where $(\zkl)_{m\geq 1}$ is a sequence of independent and identically distributed Poisson point processes on $E_0$, each having intensity measure $\t_k$.
\enpr

To prepare the proof, some elements of which are similar to the main computations in \cite{otto:thale:2022}, we recall that for $\omega\su \mathbb H^d$ locally finite, we put $f(x, \omega) := |B_{R_k(x, \omega)}|_{\ms{hyp}} - v_\la$. We also define $g(x, \omega) := \one\{f(x, \omega) > s_0\}$ and $\mathcal S(x,\omega)=B_{R_k(x,\omega)}(x)$. Then, $f$ and $g$ are localized to $\mathcal S$ in the sense of \cite{bobrowski:schulte:yogeshwaran:2021}. Namely, for every $x\in \omega$ and all $S\supseteq \mathcal S(x,\omega)$, we have that 
$
g(x,\omega) = g(x,\omega \cap S), \text{ and }f(x,\omega) = f(x,\omega \cap S) \text{ if } g(x,\omega)=1. 
$
Moreover, $\mathcal S(x,\omega)$ is a so-called stopping set, in the sense that if $B_{R_k(x,\omega)}(x) \su S$, then $B_{R_k(x,\omega \cap S)}(x) \su S$
for every compact set $S\su \mathbb H^d$. 
First, we set $S_x:= B_{r_\la (w_\la)}(x)$. Henceforth, we put 
$$s_\la :=\P \big( \PP ( B_{r_\la(s_0)}(x) ) \le k-1 \big)\quad\text{ and }\quad s_\la' :=\P \big( \PP ( B_{r_\la(w_\la)}(x)) \le k-1 \big).$$
\begin{proof}[Proof of Proposition \ref{lem:sep}]
We shall check a set of sufficient conditions for convergence in Kantorovich-Rubinstein distance as provided in Theorem 6.4 in \cite{bobrowski:schulte:yogeshwaran:2021}. This requires the analysis of the total variance distance $\dist_{\ms{TV}}(\E [\eta_{\la}^{(m)}], \tau_k)$ as well bounds on three error terms $E_1$, $E_2$ and $E_3$ which will be defined below.

First, it is claimed that 
$$
\dist_{\ms{TV}} \big(\E [\eta_{\la}^{(m)}], \tau_k\big)\to 0, \ \ \text{as } \la \uparrow\infty. 
$$
	Recall that $f(x, \PP_{Q_m}) > u$ if and only if $\PP \big( B_{r_\la(u)}(x) \big) \le k$. Hence, for $u >s_0$, it follows from the Mecke equation for Poisson point processes \cite[Theorem 4.1]{LastPenroseBook}, that 
\begin{align*}
\E \big[ \eta_{\la}^{(m)}(u,\infty) \big] &= \E \Big[ \sum_{x\in \PP_{Q^-_m}} \one \big\{f(x,\PP_{Q_m}) >u \big\} \Big] = \E\Big[ \sum_{x\in \PP_{Q^-_m}} \one \big\{\PP \big( B_{r_\la(u)}(x) \big) \le k \big\} \Big] \\
&=\int_{Q^-_m} \P \big( \PP \big( B_{r_\la(u)}(x) \big) \le k-1 \big)V_{\ms{hyp}}(\dif x) = |Q^-_m|_{\ms{hyp}} \sum_{m=0}^{k-1} e^{-(u+v_\la)} \frac{(u+v_\la)^m}{m !}. 
\end{align*}
This implies that $\E[\eta_{\la}^{(m)}]$ has the Lebesgue density 
$$
q_k(u) := |Q^-_m|_{\ms{hyp}}\, \frac{e^{-(u+v_\la) }(u+v_\la)^{k-1}}{(k-1)!}, \ \ \ u >s_0, 
$$
and thus, 
\begin{align*}
&\dist_{\ms{TV}}\big(\E [\eta_{\la}^{(m)}], \, \tau_k \big) \le \int_{s_0}^\infty \Big| q_k(u)-\frac{e^{-u}}{(k-1)!} \Big| \dif u \\
&\le \big|1 -|Q^-_m|_{\ms{hyp}} e^{-v_\la} v_\la^{k-1} \big| \int_{s_0}^\infty \Big( 1+\frac{u}{v_\la} \Big)^{k-1}e^{-u}\dif u + \int_{s_0}^\infty  \Big| \Big( 1+\frac{u}{v_\la} \Big)^{k-1}-1 \Big|e^{-u}\dif u. 
\end{align*}
	The second term above vanishes as $\la \uparrow\infty$ because of the dominated convergence theorem. For the first term, we note that Proposition \ref{pr:cvkm} gives 
$$
	|Q^-_m|_{\ms{hyp}} e^{-v_\la} v_\la^{k-1} = \f{|Q^-_m|_{\ms{hyp}}}{|Q_m|_{\ms{hyp}}}|Q_m|_{\ms{hyp}}e^{-v_\la} v_\la^{k-1} = \f{|Q^-_m|_{\ms{hyp}}}{|Q_m|_{\ms{hyp}}}\to 1, \ \ \text{ as } \la \uparrow\infty. 
$$
It is now concluded that $\dist_{\ms{TV}}\big(\E [\eta_{\la}^{(m)}], \tau_k\big)\to 0$, as $\la\to\infty$.

	To ease notation, we henceforth write $r_\la$ instead of $r_\la(s_0)$.
Subsequently, we demonstrate that 
$$
E_1 := \int_{Q^-_m} \E \Big[ g(x,\PP_{Q_m} + \de_x)\, \one \big\{\mathcal S (x,\PP_{Q_m} + \de_x) \not\su S_x\big\} \Big]V_{\ms{hyp}}(\dif x) \to 0, \ \ \ \la \uparrow\infty. 
$$
Notice that 
\begin{align*}
\one \big\{\mathcal S (x,\PP_{Q_m} + \de_x) \not\su S_x\big\} = \one \big\{(\PP_{Q_m} + \de_x) \big( B_{r_\la(w_\la)}(x) \big) \le k \big\} = \one \big\{\PP \big(B_{r_\la(w_\la)}(x) \big)\le k-1 \big\}, 
\end{align*}
where the second equality is due to the fact that $B_{r_\la(w_\la)}(x)\su Q_m$ for all $x\in Q^-_m$. Now,
\begin{align*}
	E_1 &\le \int_{Q^-_m} s_\la'\,V_{\ms{hyp}}(\dif x) \le |Q_m|_{\ms{hyp}} \sum_{m=0}^{k-1} e^{-(v_\la+w_\la)} \frac{(v_\la+w_\la)^m}{m!}\le ke^{-w_\la }\Big(1+\frac{w_\la}{v_\la}\Big)^{k-1} \to 0, \ \ \ \la \uparrow\infty,
\end{align*}
where we have used that $|Q_m|_{\ms{hyp}}=e^{v_\la} v_\la^{-(k-1)}$.
 
Next, we turn our attention to the second error term
$$
E_2 := \int_{Q^-_m}\int_{Q^-_m} \one \{S_x \cap S_z \neq \emptyset \} \E \big[ g(x, \PP_{Q_m} + \de_x) \big]\E \big[ g(z, \PP_{Q_m} + \de_z) \big]V_{\ms{hyp}}(\dif z)V_{\ms{hyp}}(\dif x).
$$
We here apply the following inequalities: 
	\begin{align*}
		\one \{S_x\cap S_z\neq \emptyset \}&\le \one \big\{\dist_{\ms{hyp}}(x,z) \le 2 r_\la (w_\la) \big\}, \\
		\E \big[ g(x, \PP_{Q_m} + \de_x) \big] &= s_\la = \sum_{i=0}^{k-1} \frac{e^{-(s_0+v_\la)} (s_0+v_\la)^i}{i!} \le Cv_\la^{k-1}e^{-v_\la}. 
	\end{align*}
Hence, 
\begin{align}
E_2 &\le C(v_\la^{k-1}e^{-v_\la})^2 \int_{Q^-_m}\int_{\mathbb H^d}\one \big\{\dist_{\ms{hyp}}(x,z) \le 2 r_\la (w_\la) \big\}V_{\ms{hyp}}(\dif z)V_{\ms{hyp}}(\dif x) \label{e:E2}\\
&= C(v_\la^{k-1}e^{-v_\la})^2 |Q^-_m|_{\ms{hyp}} | B_{2r_\la(w_\la)} |_{\ms{hyp}} \notag\\
&\le C(v_\la^{k-1}e^{-v_\la})^2 e^{v_\la} v_\la^{-(k-1)}\, e^{2r_\la(w_\la)(d-1)}\notag \\
&\le C v_\la^{k-1}v_\la^2 e^{-v_\la} \to 0, \ \ \text{as } \la \uparrow\infty. \notag
\end{align}
For the inequality at the last line, we have applied \eqref{e:bound.r.la.w.la}.

For the bound of $E_3$ we distinguish between the cases $k=1$ and $k \ge2$. If $k=1$ we have that
$$\E \big[ g(x, \PP_{Q_m} + \de_x + \de_z)g(z, \PP_{Q_m} + \de_x + \de_z) \big] \le s_\lambda^{3/2}$$
and, hence, 
\begin{align*}
	E_3&:= \int_{Q^-_m}\int_{Q^-_m} \one \{S_x \cap S_z \neq \emptyset \} \E \big[ g(x, \PP_{Q_m} + \de_x + \de_z)g(z, \PP_{Q_m} + \de_x + \de_z) \big]V_{\ms{hyp}}(\dif z)V_{\ms{hyp}}(\dif x)\\
	&\le s_\lambda^{3/2} \int_{Q^-_m}\int_{\mathbb H^d}\one \big\{\dist_{\ms{hyp}}(x,z) \le 2 r_\la (w_\la) \big\}V_{\ms{hyp}}(\dif z)V_{\ms{hyp}}(\dif x)\\
	&\le C(e^{-v_\la})^{3/2} \frac{e^{\la (d-1)}}{u_\la}\, e^{2r_\la(w_\la)(d-1)}\le C v_\la^2 e^{-v_\la/2} \to 0, \ \ \text{as } \la \uparrow\infty,
\end{align*}
where in the last step we used the volume estimate for $Q_m^-$ and the fact that $|B_r(x)|_{\ms{hyp}}\leq Ce^{(d-1)r}$ for any $r>0$ and $x\in\H^d$; see the discussion after \eqref{eq:VolBall}.

 For $k \ge 2$, we split $E_3$ into two terms: 
\begin{align*}
E_3 &= \int_{Q^-_m}\int_{Q^-_m} \one \{S_x \cap S_z \neq \emptyset \} \E \big[ g(x, \PP_{Q_m} + \de_x + \de_z)g(z, \PP_{Q_m} + \de_x + \de_z) \big]V_{\ms{hyp}}(\dif z)V_{\ms{hyp}}(\dif x)\\
&\le \int_{Q^-_m}\int_{Q^-_m} \one \big\{\dist_{\ms{hyp}}(x,z) \le ar_\la(w_\la) \big\} \\
&\qquad \times \P\Big( (\PP_{Q_m} + \de_z)\big( B_{r_\la} (x)\big) \le k-1,\, (\PP_{Q_m} + \de_x)\big( B_{r_\la} (z)\big) \le k-1 \Big) V_{\ms{hyp}}(\dif z)V_{\ms{hyp}}(\dif x) \\
&+ \int_{Q^-_m}\int_{Q^-_m} \one \big\{ar_\la(w_\la)< \dist_{\ms{hyp}}(x,z) \le 2r_\la(w_\la) \big\} \\
&\qquad \times \P\Big( (\PP_{Q_m} + \de_z)\big( B_{r_\la} (x)\big) \le k-1,\, (\PP_{Q_m} + \de_x)\big( B_{r_\la} (z)\big) \le k-1 \Big) V_{\ms{hyp}}(\dif z)V_{\ms{hyp}}(\dif x) \\
&=: E_{3,1}+E_{3,2},
\end{align*}
where $0<a<1$ is so small that $2r_\la-a(k-1)r_\la(w_\la)>0$. The term $E_{3,1}$ is bounded by 
\begin{align*}
E_{3,1} &\le \int_{Q^-_m}\int_{Q^-_m} \one \big\{\dist_{\ms{hyp}}(x,z) \le ar_\la(w_\la) \big\}\P\big( \PP \big( B_{r_\la}(x) \big)\le k-2 \big)\\
&\qquad \qquad \qquad\times \P\big( \PP \big( B_{r_\la}(z)\setminus B_{r_\la}(x) \big)\le k-2 \big)V_{\ms{hyp}}(\dif z)V_{\ms{hyp}}(\dif x). 
\end{align*}
Let $p \in \H^d$ be a fixed reference point (sometimes called the origin) and let $\exp_p:T_p\to\H^d$ denote the exponential map at $p$ in which $T_p\cong\R^d$ is the tangent space at $p$. Then, by the above
the bound,
$$
 \P\big( \PP \big( B_{r_\la}(x) \big)\le k-2 \big)\le Ce^{-v_\la}v_\la^{k-2}.
$$
 It follows from the polar integration formula in hyperbolic geometry \cite[pp. 123-125]{flavors} that
\begin{align*}
E_{3,1} &\le Ce^{-v_\la} v_\la^{k-2} |Q_m|_{\ms{hyp}} \sum_{i=0}^{k-2} \int_{\mathbb H^d} \one \big\{\dist_{\ms{hyp}}(z,p) \le ar_\la(w_\la) \big\} \\
&\qquad \qquad \qquad\qquad \qquad \times e^{-| B_{r_\la}(z)\setminus B_{r_\la}(p)|_{\ms{hyp}}} \frac{\big( |B_{r_\la}(z)\setminus B_{r_\la}(p)|_{\ms{hyp}}\big)^i}{i!}V_{\ms{hyp}}(\dif z)\\
&\le \frac{C}{v_\la} \sum_{i=0}^{k-2}\int_0^{ar_\la(w_\la)} \sinh^{d-1}(s) e^{-| B_{r_\la}(z)\setminus B_{r_\la}(p)|_{\ms{hyp}}} \frac{\big( |B_{r_\la}(z)\setminus B_{r_\la}(p) |_{\ms{hyp}}\big)^i}{i!}\dif s, 
\end{align*}
where $z=\exp_p(sv_0)$ for some $v_0\in\mathbb{S}_p^{d-1}$, the $(d-1)$-dimensional unit sphere in $T_p$. It follows from Lemma 5 in \cite{otto:thale:2022} that 
\begin{align*}
E_{3,1} &\le \frac{C}{v_\la} \sum_{i=0}^{k-2}\int_0^{ar_\la(w_\la)} \Big( \frac{e^s}2\Big)^{d-1} e^{-\alpha_1 s e^{(d-1)(r_\la-s/2)}} \frac{\big( \alpha_2 s e^{(d-1)r_\la}\big)^i}{i!}\dif s\\
&\le \frac{C}{v_\la} \sum_{i=0}^{k-2} \frac{e^{(d-1)ir_\la}}{i!}\, \int_0^\infty e^{-s(\alpha_1e^{(d-1)(2r_\la-ar_\la(w_\la))/2}-(d-1))}\,s^i \dif s \\
&\le \frac{C}{v_\la} \sum_{i=0}^{k-2} e^{-(d-1)\frac{2r_\la-a(i+1)r_\la(w_\la)}2} \\
&\le \frac{C}{v_\la} e^{-(d-1)\frac{2r_\la-a(k-1)r_\la(w_\la)}2} \le \frac{C}{v_\la}\to 0, \ \ \text{as } \la \uparrow\infty. 
\end{align*}
Here, the last inequality follows from the constraint $2r_\la-a(k-1)r_\la(w_\la)>0$. 

For $E_{3,2}$, one can see that 
\begin{align*}
E_{3,2} &\le \int_{Q^-_m}\int_{Q^-_m} \one \big\{ar_\la(w_\la)< \dist_{\ms{hyp}}(x,z) \le 2r_\la(w_\la) \big\} s_\la \\
&\qquad\times\P \big(\PP \big( B_{r_\la} (z) \setminus B_{r_\la} (x)\big)\le k-1 \big)V_{\ms{hyp}}(\dif z)V_{\ms{hyp}}(\dif x) \\
&\le Cs_\la|Q_m|_{\ms{hyp}} \int_{\mathbb H^d} \one \big\{ar_\la(w_\la)< \dist_{\ms{hyp}}(z,p) \le 2r_\la(w_\la) \big\} \\
&\qquad\times
 \sum_{i=0}^{k-1} e^{-| B_{r_\la} (z) \setminus B_{r_\la}(p) |_{\ms{hyp}}} \frac{ | B_{r_\la} (z) \setminus B_{r_\la}(p) |_{\ms{hyp}} ^i}{i!} V_{\ms{hyp}}(\dif z)\\
&\le C\sum_{i=0}^{k-1}\frac1{i!}\,\int_{ar_\la(w_\la)}^{2r_\la(w_\la)} \sinh^{d-1}(s)\, e^{-| B_{r_\la} (z) \setminus B_{r_\la}(p) |_{\ms{hyp}}} \big( | B_{r_\la} (z) \setminus B_{r_\la}(p) |_{\ms{hyp}} \big)^i \dif s, 
\end{align*}
where again $z=\exp_p(sv_0)$ for some $v_0\in\mathbb{S}_p^{d-1}$. 
Note that for every $z=\exp_p(sv_0)\in \mathbb H^d$ with $ar_\la(w_\la)\le s \le 2r_\la(w_\la)$, 
$$
| B_{r_\la} (z) \setminus B_{r_\la}(p) |_{\ms{hyp}} \ge |B_{s/2}(p)|_{\ms{hyp}}, 
$$
and further, by Lemma 4 in \cite{otto:thale:2022}, $|B_{s/2}(p)|_{\ms{hyp}}\ge C_1 e^{s(d-1)/2}$.
Using this bound, 
\begin{align*}
E_{3,2} &\le C\sum_{i=0}^{k-1}\frac1{i!}\,\int_{ar_\la(w_\lambda)}^{2r_\la(w_\lambda)} e^{s(d-1)} e^{- C_1 e^{s(d-1)/2}} \big( C_1 e^{s(d-1)/2} \big)^i \dif s\\
&\le Ce^{-(d-1)ar_\la(w_\la)/2} \sum_{i=0}^{k-1}\frac{C_1^{i}}{i!}\,\int_{ar_\la(w_\lambda)}^{2r_\la(w_\lambda)} e^{- C_1 e^{s(d-1)/2}} (e^{s(d-1)/2})^{i+3} \dif s.
\end{align*}
Applying the substitution $u=e^{s(d-1)/2}$, we find that the integral is bounded in $\lambda$ and, hence, that $E_{3,2}\to 0$, as $\lambda \to \infty$.

We now put all bounds together and apply Theorem 6.4 in \cite{bobrowski:schulte:yogeshwaran:2021} to conclude that, for each $m\geq 1$,
$$
\dist_{\ms{KR}}\big(\LL(\ekl), \LL(\zkl)\big) \leq \dist_{\ms{TV}} \big(\E [\eta_{\la}^{(m)}], \tau_k\big) + 2(E_1+E_2+E_3)\to 0,, \ \ \text{as } \la \uparrow\infty.
$$
This completes the proof.
\enp

%
%
By the maximal coupling lemma \cite[Lemma 4.32]{randMeas}, Proposition \ref{lem:sep} gives for each $m\geq 1$ couplings $\hekl$ and $\hzkl$ satisfying 
\begin{align}
	\label{eq:coup}\P\big(\hekl \ne \hzkl\big) = \dist_{\ms{TV}}\big(\LL(\ekl), \LL(\zkl)\big) \to 0, \ \ \text{as } \la \uparrow\infty.
\end{align}
Hence, we can define a coupling between $\ek$ and $\zk$ by setting 
$$ \hek := \sum_{m \le u_\la}\hekl, \quad\text{ and }\quad\hzk := \sum_{m \le u_\la}\hzkl, $$
where we assume that the pairs $(\hekl, \hzkl)$ are independent and identically distributed for different $m$.  Having constructed the couplings, we now proceed with the proof of the exponential equivalence.

%
%
\bep[Proof of Proposition \ref{pr:ee1}]
First, we note that for every $a > 0$ 
\begin{align*}
\f1{u_\la}\log\P\big(\dist_{\ms{TV}}\big(\hek, \hzk\big) \ge \de u_\la\big) &\le \f1{u_\la}\log\P\Big(\sum_{m \le u_\la}\dist_{\ms{TV}}\big(\hekl, \hzkl\big) \ge \de u_\la\Big) \\
&\leq{1\over u_\la}\log\Big(e^{-a\delta u_\la}\E\Big[ e^{a\sum_{m \le u_\la}\dist_{\ms{TV}}\big(\hekl, \hzkl\big) }\Big]\Big)\\
&= -a\delta +{1\over u_\la}\sum_{m\leq u_\la}\log \E[e^{a\,\dist_{\ms{TV}}\big(\hekl, \hzkl\big)}]\\
&\le -a\de + \log\E[e^{a\,\dist_{\ms{TV}}(\heol, \hzol)}],
\end{align*}
where we used the exponential Markov inequality and the i.i.d.~property of the pairs $(\hekl, \hzkl)$ discussed above.
Moreover, by \eqref{eq:coup}, we have that $\dist_{\ms{TV}}(\heol, \hzol) \to 0$, as $\la\uparrow\infty$ in probability. Hence, as in \cite[Proposition 4.3]{kang}, it suffices to verify the uniform integrability condition 
\begin{align}
	\label{eq:ee11}
	\limsup_{\la\tff}\E\Big[e^{a\QQ(Q^-_1)}\Big]< \ff
\end{align}
for every $a > 0$, where we let
$$\QQ := \sum_{X_i\in\Pn   }\one\{\Pn (B_{\rl}(X_i)) \le k\}\de_{X_i}$$
denote the point process of exceedances. The key observation that we will use tacitly throughout the rest of the proof is that to determine the restriction of $\QQ$ to a set $A\su \H^d$, it suffices to know $\Pn$ in an $\rl$-neighborhood of $A$. Therefore, we partition $Q_1$ into vertical layers $Q_{1, \ell} := S_1\ti T_\ell$, where $T_\ell := [e^{-\la + \ell \rl}, e^{-\la + (\ell + 1)\rl})$. In order to avoid dealing with an infinite number of such layers, we define a value $\ell_0\ge 1$ such that $S_1 = [0, e^{-\la +  (\ell_0 + 2) \rl}]^{d - 1}$. As in previous arguments, here, $\ell_0$ is assumed to be an integer; if this is not the case, we replace $\rl$ by $r'$ which is a small perturbation of $\rl$.  We then set $Q_{1, \ff} := S_1 \ti [e^{-\la + (\ell_0+1) \rl},\ff)$.
In particular, we deduce from Lemma \ref{lem:dist} that $\QQ_{Q_{1, \ell}}$ is independent of $\QQ _{Q_{1, \ell'}}$ whenever $\ell$ and $\ell'$ are such that $|\ell' - \ell|\ge 3$. Hence, it suffices to prove that 
\begin{align}
	\label{eq:uint}
	\limsup_{\la\tff}  \E \Big[ e^{a\QQ(Q_{1, \ff}) }\Big] < \ff \quad \text{ and }\quad \limsup_{\la\tff} \prod_{0 \le \ell\le \ell_0, \ell \in 3\Z} \E \Big[ e^{a\QQ(Q_{1, \ell}) }\Big] < \ff.
\end{align}
The arguments for layers with $\ell -1 \in 3\Z$ or $\ell - 2 \in 3\Z$ are similar. 

%
%
We start with $Q_{1, \ff}$, where to ease notation, we set $J := Q_{1, \ff}$. Note that by the kNN property, we have that if $X_i$ is not among the $(k - 1)$ nearest neighbors of $X_j$, then 
$$\big|B(X_i, X_j)\big|_{\ms{hyp}}:= \big|B_{\rl/2}(X_i)\sm B_{\rl}(X_j)\big|_{\ms{hyp}}\ge {1\over 2}|B_{\rl/2}(X_i)|_{\ms{hyp}}.$$
By the same reason, each point $z \in B_{\rl/2}(X_i)$ is covered by at most $(k - 1)$ further balls of the form $B_{\rl/2}(X_{i'})$. Therefore,  
\begin{align*}
	\Big|\bigcup_{X_i \in \QQ^{}_{J\sm B_{\rl}(X_1)}} \hspace{-0.2cm}B(X_i, X_1)\Big|_{\ms{hyp}}
	\ge  {1\over{2 k}}(\QQ(J) - k) |B_{\rl/2}(X_i)|_{\ms{hyp}}
	\ge   C_0(\QQ(J) - k)  e^{(d - 1)\rl/2},
\end{align*}
 where in the last step we used once again Lemma 4 in \cite{otto:thale:2022}. Now, note that $|J|_{\ms{hyp}} \in O\big(e^{(d - 1)\rl}\big)$ and  let for $n\in\mathbb{N}$, $Z_n$ be a Poisson variable with mean $C_0(n - k)_+  e^{(d - 1)\rl/2}$.  Thus, by the multivariate Mecke equation for Poisson point processes \cite[Theorem 4.4]{LastPenroseBook}, for  $n\in\mathbb{N}$,
\begin{align}
	\notag\P(\QQ(J) = n) &\le \f1{n!} \E\Big[\hspace{-.6cm}\sum_{\substack{X_1, \dots, X_n \in \Pn\cap J\\ \text{ pairwise distinct}}}\hspace{-.6cm}\one\{\Pn(B_{\rl}(X_1)) \le k \}\one\Big\{\Pn\Big(\hspace{-.5cm}\bigcup_{X_i \not \in B_{\rl}(X_1)} B(X_i, X_1)\Big)\le  kn\Big\}\Big]\\
	&\le \f1{n!} s_\la |J|_{\ms{hyp}}^n\P\big(Z_n\le  kn\big).\label{eq:1423A}
\end{align}
Now, we want to apply the concentration inequality for Poisson random variables \cite[Lemma 1.2]{penrose}. That means, we need to provide an upper bound for the expression 
$$\exp\Big(- kn \log\Big(\f{kn}{eC_0(n - k)  e^{(d - 1)\rl/2}}\Big)\Big)\exp\Big( - C_0(n - k)  e^{(d - 1)\rl/2}\Big).$$
To that end, we note that the first exponential is bounded by 
$\big( C_1e^{ (d - 1)\rl/2}\big)^{km}$
 for a suitable $C_1 > 0$. Moreover, since $\log(|J|_{\ms{hyp}}) \in O(\rl)$, we conclude that $\log(|J|_{\ms{hyp}}) \in o\big(e^{(d - 1)\rl/2}\big)$.
 Plugging this into \eqref{eq:1423A} gives that for $n \ge 2k$ and suitable $C_2 >0$, 
 \begin{align*}
	\P(\QQ(J) = n) 
	&\le \f{|J|_{\ms{hyp}}}{n!} s_\la\big(C_1 e^{ (d - 1)\rl/2}\big)^{kn}\exp\Big(-C_0(n - k) e^{(d - 1)\rl/2}\Big)\\
	&\le \f{|J|_{\ms{hyp}}}{n!}s_\la \exp(-C_2n e^{(d - 1)\rl/2}).
\end{align*}
Therefore, by Markov's inequality,
\begin{align*}
		 \E \big[ e^{a\QQ(J) }\big] 
	\le  1 + e^{2ka}|J|_{\ms{hyp}}s_\la + \sum_{m \ge 2k + 1}\f{|J|_{\ms{hyp}}}{m!}s_\la e^{-C_2m e^{(d - 1)\rl/2}}
	\le 1 + ( e^{2ka} + 1)|J|_{\ms{hyp}}s_\la.
\end{align*}
Hence, for large $\la$ we have 
$$\E \big[ e^{a\QQ(J) }\big]  \le e^{( e^{2ka} + 1) |J|_{\ms{hyp}}s_\la},$$
which remains bounded since $|J|_{\ms{hyp}}s_\la \in O(1)$.

%
%
Next, to bound $\E \big[ e^{a\QQ(Q_{1, \ell})}\big]$, we proceed similarly as in the proof of \cite[Proposition 4.3]{kang} and decompose $Q_{1, 3\ell}$ into independent horizontal blocks. More precisely, we consider the diluted family of boxes $G_\ell := \big\{8e^{-\la + (\ell + 1) \rl}z + J_\ell \su Q_{1, 3\ell}\big\}$, where we set 
$$J_\ell := [0,e^{-\la + (\ell + 1)\rl}]^{d - 1} \ti T_\ell.$$
Hence, applying H\"older's inequality we reduce the verification of \eqref{eq:uint} to showing that
\begin{align}
	\label{eq:uint1}
	\limsup_{\la\tff} \prod_{\ell \le \ell_0} \E \big[ e^{a\QQ(J_\ell) }\big]^{\dll} < \ff,
\end{align}
where $\dll :=  |Q_{1, \ell}|_{\ms{hyp}}/|J_{\la, \ell}|_{\ms{hyp}}$. Now, proceeding as in the case of $Q_{1, \ff}$ shows that for large $\la$ 
\begin{align*}
	\prod_{\ell \le \ell_0} \E \big[ e^{a\QQ(J_\ell) }\big]^{\dll}\le e^{( e^{2ka} + 1)\dll |J_{\la, \ell}|_{\ms{hyp}}s_\la}\le e^{( e^{2ka} + 1) |Q_1|_{\ms{hyp}}s_\la}.
\end{align*}
Since the last expression belongs to $O(1)$, the proof is complete.
\enp

%
%
\section{Proof of Proposition \ref{pr:ee2}}
\label{sec:ee2}
In this section, we prove Proposition \ref{pr:ee2}, that is,  the exponential equivalence of the families of random measures $\ek/u_\la$ and $\xi_{ \la}/u_\la$. To that end, we will  show that two error terms are negligible. More precisely, we set 
$$\nla := \#\big\{x \in \Pn \cap W_\la^- \co R_k(x) > \rlw\big\}$$
and 
$$\nlb := \#\big\{x \in \Pn \cap (W_\la \sm W_\la^-) \co R_k(x) > \rl\big\}.$$
where $W_\la^- := \bigcup_{m\leq u_\la} (Q_m \sm Q^-_m)$.
Then, we show that $\nla$ and $\nlb$ are exponentially negligible.

%
%
\bel[$\nla$ is exponentially negligible]
\label{lem:ee2a}
The family of random variables $\nla/u_\la$ is exponentially equivalent to $0$.
\enl

%
%
\bel[$\nlb$ is exponentially negligible]
\label{lem:ee2b}
The family of random variables $\nlb/u_\la$ is exponentially equivalent to $0$.
\enl

%
%
Before proving Lemmas \ref{lem:ee2a} and \ref{lem:ee2b}, we explain how to conclude the proof of Proposition \ref{pr:ee2}.

\bep[Proof of Proposition \ref{pr:ee2}]
First, recall that 
$$\ek := \sum_{m \le u_\la}\sum_{x \in \PP_{Q^-_m}}\one\{f(x, \PP_{Q_m}) > s_0\}\delta_{f(x, \PP_{Q_m})}. $$
Now, since $x\in Q^-_m$, Lemma \ref{lem:ee2a} gives that $\ek$ is exponentially equivalent to 
$$\ek' := \sum_{m \le u_\la}\sum_{x \in \PP_{Q^-_m}}\one\{f(x, \Pn) > s_0\}\delta_{f(x, \Pn)}.$$
Moreover, Lemma \ref{lem:ee2b} gives that $\ek'/u_\la$ is exponentially equivalent to
$$\sum_{x \in \Pn_{W_\la}}\one\{f(x, \Pn) > s_0\}\delta_{f(x, \Pn)}.$$
Noting that the latter expression coincides with $\xi_{ \la}$ concludes the proof.
\enp

%
%
It remains to establish Lemmas \ref{lem:ee2a} and \ref{lem:ee2b}. We begin with the proof of Lemma \ref{lem:ee2a}, which is similar to that of \cite[Proposition 4.4]{kang}. Here, we stress that similarly as in the proof of \cite[Proposition 4.4]{kang}, the families
$\big\{x \in \Pn_{Q^-_m}\co R_k(x) > \rlw\big\}$
are independent for different $m\geq 1$.

\bep[Proof of Lemma \ref{lem:ee2a}]
Recalling the notation $\dll = |Q_{1, \ell}|_{\ms{hyp}}/|J_{\ell}|_{\ms{hyp}}$, we proceed as in the proof of Proposition \ref{pr:ee1}. More precisely, we need to show that for all $a > 0$, 
\begin{align}
        \label{eq:uint2}
        \limsup_{\la\tff} \prod_{m \ge 0} \E \big[ e^{a\QQ'(J_\ell) }\big]^{\dll} \le 1,
\end{align}
where $J_\ell := [0,e^{-\la + (\ell + 1)\rl}]^{d - 1} \ti T_\ell$ and 
$$\QQ' := \sum_{X_i\in\Pn } \one \big\{ \Pn (B_{\rlw}(X_i)) \le k \big\}\de_{X_i}.$$
By copying the arguments from the proof of Proposition \ref{pr:ee1}, we arrive at
        $$\limsup_{\la\tff} \prod_{m \ge 0} \E \big[ e^{a\QQ'(J_\ell) }\big]^{\dll} \le e^{( e^{2ka} + 1) |Q_1|_{\ms{hyp}}s_\la'}.$$
Hence, noting that $ |Q_1|_{\ms{hyp}}s_\la'\in o(1)$ concludes the proof of \eqref{eq:uint2}.
\enp

%
%
Finally, we prove Lemma \ref{lem:ee2b}. Here, we will use a diluted family of cubes. 
\bep[Proof of Lemma \ref{lem:ee2b}]
The idea is to proceed similarly as in the proof of Proposition \ref{pr:ee1}.   Indeed, by Markov's inequality, it suffices to show that for every $a > 0$ we have that 
$$\log\E \big[ e^{a\QQ(W_\la^-) }\big] \in o(u_\la).$$
As in the proof of Proposition \ref{pr:ee1}, we subdivide $W_\la$ into vertical layers $W_{\la, \ell} := [0, 1]^{d - 1}\ti T_\ell$, where $T_\ell := [e^{-\la + \ell \rl}, e^{-\la + (\ell + 1)\rl}]$. We also set $W_{\la, \ell}^- := W_{\la, \ell} \cap W_\la^-$.
In particular, we deduce from Lemma \ref{lem:dist} that $\QQ \cap W_{\la, \ell}^-$ is independent of $\QQ \cap W_{\la, \ell'}'$ for all $\ell$ and $\ell'$ satisfying $|\ell' - \ell|\ge 3$. Again, we define a value $\ell_0\ge 1$ such that $[0, 1]^{d - 1} = \big[0, e^{-\la +  (\ell_0 + 2) \rl}\big]^{d - 1}$.
Hence, it suffices to prove that for every $a > 0$,
\begin{align}
	\label{eq:uint2a0}
	 \sum_{\ell \ge \ell_0 + 1, \ell \in 3\Z}\log\E \big[ e^{a\QQ(W_{\la, \ell}^-) }\big] \in o(u_\la)\quad \text{ and }\quad   \sum_{0 \le \ell\le \ell_0, \ell \in 3\Z} \log\E \big[ e^{a\QQ(W_{\la, \ell}^-) }\big] \in o(u_\la).
\end{align}
The arguments for layers with $\ell - 1\in 3\Z$ or $\ell - 2 \in 3\Z$ are similar. The key step in the proof of \eqref{eq:uint2a0} is the exponential moment bound
\begin{align}
	\label{eq:uint2a1}
	 \log\E \big[ e^{a\QQ(W_{\la, \ell}^-) }\big] \in O\big(u_\la| W_{\la, \ell}^-|_{\ms{hyp}}/|W_\la|_{\ms{hyp}}\big).
\end{align}
Once \eqref{eq:uint2a1} is shown, \eqref{eq:uint2a0} follows by summation over $\ell \ge 0$.  We first recall from the proof of Proposition \ref{pr:ee1} that 
\begin{align}
	\label{eq:emo}
	\E \big[ e^{a\QQ(W_\la^- \cap I) }\big] \le  1 + ( e^{2ka} + 1)|W_\la^- \cap I|_{\ms{hyp}}s_\la,
\end{align}
where $I$ is a set of the form $I = S \ti T_\ell$ with $|I|_{\ms{hyp}}\in O\big(e^{(d - 1)\rl}\big)$ for some $S \su \R^{d - 1}$.

%
%
Now, to prove \eqref{eq:uint2a1}, we first consider the case where $\ell \ge \ell_0 + 1$. Then, by \eqref{eq:emo}, 
$$\log\E \big[ e^{a\QQ(W_{\la, \ell}^-) }\big] \le{( e^{2ka} + 1)} | W_{\la, \ell}^-|_{\ms{hyp}}s_\la=  ( e^{2ka} + 1)u_\la| W_{\la, \ell}^-|_{\ms{hyp}}/|W_\la|_{\ms{hyp}},
$$
as asserted.

%
%
It remains to deal with the case where $\ell \le \ell_0$. Then, we proceed as in the proof of Proposition \ref{pr:ee1} and consider the diluted family of cubes $G_\ell := \big\{8e^{-\la + (\ell + 1) \rl}z + J_\ell \su W_{\la, \ell}\big\}$, where 
$$J_\ell := [0,e^{-\la + (\ell + 1)\rl}]^{d - 1} \ti T_\ell.$$
Now, we conclude from \eqref{eq:emo} that
\begin{align}
	\label{eq:uint2a2}
	 \log\E \big[ e^{a\QQ(W_\la^-\cap (J_\ell + x)) }\big] \in O\big(u_\la|W_\la^- \cap (J_\ell + x)  |_{\ms{hyp}}/|W_\la|_{\ms{hyp}}\big)
\end{align}
for any $x \in \R^{d - 1}$.  Hence, by H\"older's inequality and independence, we obtain that 
$$\log\E \big[ e^{a\QQ(W_{\la, \ell}^-) }\big] \le C_1\sum_{J \in G_\ell}\log\E \big[ e^{C_2a\QQ(W_\la^- \cap J) }\big].$$
Since any two distinct $J_1, J_2 \in G_\ell$ are disjoint, we obtain that 
$$\sum_{J \in G_\ell}|W_\la^- \cap J  |_{\ms{hyp}} = \Big|W_\la^- \cap \bigcup_{J \in G_\ell} J  \Big|_{\ms{hyp}} \le |W_{\la,\ell}'|_{\ms{hyp}} . $$
Therefore, we get 
	$\log\E \big[ e^{a\QQ(W_{\la, \ell}^-) }\big] \in O\big(|W_{\la, \ell}^-|_{\ms{hyp}}u_\la/|W_\la|_{\ms{hyp}}\big),$
	 as asserted.
\enp

\subsection*{Acknowledgement}

CT was supported by the DFG priority program SPP 2265 \textit{Random Geometric Systems}.

\bibliography{./lit}

\end{document}